\documentclass[12pt]{article}

\usepackage[utf8]{inputenc}
\usepackage[english,russian]{babel}
\usepackage{amsmath,amssymb,latexsym,amsfonts}

\begin{document} 

\begin{center}
{\bf О лакунарности и единственности для p-ичных аналогов хаоса Радемахера}
\end{center}

\begin{center}
А.~Д.~Казакова, М.~Г.~Плотников
\end{center}

\begin{center}
anna.kazakova@math.msu.ru, mikhail.plotnikov@math.msu.ru 
\end{center}
\begin{abstract}
Для двух систем функций, являющихся p-ичными аналогами хаосов Радемахера, 
доказывается их q-лакунарность и изучаются множества единственности.

Ключевые слова: q-лакунарность, неравенство Хинчина, 
множества единственности, хаос Радемахера,  
системы Виленкина--Крестенсона.

For two systems of functions that are p-ary analogues of Rademacher chaos, their q-lacunarity is proved and their uniqueness sets are studied.

\end{abstract}

\section*{Введение} 

В работе изучаются вопросы, связанные с лакунарностью 
и единственностью для систем функций. 
Систему функций $( \varphi_k \colon [ 0, 1) \to \mathbb{C}, 
\; k \in \mathbb{N} )$
назовем {\it системой единственности}, 
если для некоторого $\varepsilon > 0$ 
сходимость к нулю на множестве $E$ с 
$\mu(E) > 1 - \varepsilon$ ряда по системе $( \varphi_k )$ 
влечет равенство нулю всех его коэффициентов. 
Для фиксированного $\varepsilon > 0$ 
будем говорить, следуя~[1],  
что $( \varphi_k )$   
является {\it системой $\varepsilon$-единственности}, 
если она удовлетворяет условию из определения системы единственности 
с заданным $\varepsilon$.

Классические полные в $L_2$ системы 
функций (тригонометрическая система, системы Уолша, Хаара, Франклина) 
не являются системами единственности~[2--7]. 
Для этих систем сходимость к нулю даже почти 
всюду не гарантирует равенства нулю всех их коэффициентов 
(для тригонометрической системы примером является нуль-ряд Меньшова). 
Системами единственности 
обычно являются лакунарные в некотором смысле системы. 
В узком смысле под лакунарными понимают 
разреженные подсистемы других систем функций. 
В широком смысле лакунарность означает~[8], 
что система функций обладает некоторыми свойствами, 
присущими системам независимых (в вероятностном смысле) функций. 

Гапошкин исследовал~[8] различные вопросы, в том числе 
связанные с $\varepsilon$-единственностью, 
для систем функций, лакунарных в разных смыслах. 
В~[8] рассматривались как 
общие системы функций, так и конкретные, 
в частности, подсистемы тригонометрической системы 
и системы Уолша. 
Система функций $( \varphi_k \colon [ 0, 1) \to \mathbb{C}, 
\; k \in \mathbb{N} )$ 
называется {\it системой $q$-лакунарности}, где $q > 2$,  
если имеет место $L_2$-$L_q$-неравенство Хинчина  
$$
\left\|  
    \sum\limits_{k=1}^n c_k \varphi_k  
\right\|_{L_q}
\le 
\kappa  
\big\|  
    ( c_k )_{k=1}^n  
\big\|_{l_2}, 
\quad 
\text{$\kappa = \kappa(q) > 0$ не зависит от $c_k \in \mathbb{C}$}.            
\eqno{(1)}
$$ 
Можно писать $\Big\| \sum\limits_{k=1}^n c_k \varphi_k \Big\|_2$ 
вместо нормы справа, если $( \varphi_k )$ является ортонормированной 
системой или системой Рисса; последнее означает, что 
найдутся $R > r > 0$ такие, что 
$$ 
r
\big\|  
    ( c_k )_{k=1}^n  
\big\|_{l_2}
\le 
\left\|  
    \sum\limits_{k=1}^n c_k \varphi_k  
\right\|_{L_2}
\le 
R
\big\|  
    ( c_k )_{k=1}^n  
\big\|_{l_2}, 
\quad 
c_k \in \mathbb{C}. 
$$

Пусть система $( \varphi_k )_{ k \in \mathbb{N} }$ 
ортонормирована или является системой Рисса. 
Несложно показать 
(см., напр., [9, Ch.~2] для случая $q=4$ 
и системы Радемахера), 
что тогда из~(1) вытекает $L_1$-$L_2$-неравенство Хинчина
$$
\left\|  
    \sum\limits_{k=1}^n c_k \varphi_k  
\right\|_1
\ge 
\kappa_0 
\big\|  
    ( c_k )_{k=1}^n  
\big\|_{l_2}, 
\quad 
\text{$\kappa_0 = \kappa_0 (q) > 0$ не зависит от $c_k \in \mathbb{C}$}.
$$ 
Из этого, в свою очередь, следует, 
согласно результатам~[8], 
что $( \varphi_k )$ является системой единственности. 
В частности, любая система $q$-лакунарности, 
являющаяся ортонормированной или системой Рисса, 
является системой единственности. 

Классическим примером лакунарной системы 
является система функций Радемахера~[10-11]. 
С одной стороны, эта система 
является лакунарной по Адамару подсистемой 
системы Уолша, одной из классических систем 
функций в анализе.  
С другой стороны система Радемахера 
является системой независимых симметричных бернуллиевских 
случайных величин. 

Стечкин и Ульянов установили~[12], 
что система Радемахера является системой $1/2$-единственности, 
причем константа $1/2$ неулучшаема.  
Кури обобщил~[13] этот результат 
на случай произвольных лакунарных по Адамару подсистем системы Уолша. 
Проблемам единственности для системы Радемахера 
и других подсистем системы Уолша 
посвящены также работы Лукомского~[14,\,15]. 

В работах~[16,\,17] 
Бонами доказала  (см. также~[9,\,18]), 
что система, состоящая из 
$d$-членных произведений различных функции Радемахера 
({\it $d$-хаос Радемахера}) является системой $q$-лакунарности. 
Значит, она является и системой единственности. 
Асташкин и Суханов нашли~[1], 
что точная константа $\varepsilon$-единственности 
для $d$-хаоса Радемахера есть $\varepsilon = 1/2^d$. 

В нашей работе мы обобщаем результаты о 
$q$-лакунарности и $\varepsilon$-единственности 
со случая двоичной системы счисления, 
с которой связана система Радемахера, 
на случай $p$-ичной системы с произвольным 
натуральным основанием $p \ge 2$. 
При этом возникают как минимум два способа обобщить 
хаос Радемахера на $p$-ичный случай. 
Рассмотрим системы функций 
$$
\{ VC_n, \; n \in V^{(d)}_p \},    
\eqno{(2)}
$$
$$ 
\{ VC_n, \; n \in \widetilde{V}^{(d)}_p \}.   
\eqno{(3)}
$$
Здесь $VC_n$ --- функции Виленкина~---Крестенсона, 
$p \ge 2$ и $d \ge 1$ --- натуральные числа, 
а $V^{(d)}_p \subset \widetilde{V}^{(d)}_p$ --- множества, 
состоящие из всех натуральных $n$ соответственно 
следующего вида (все $k_j \in \mathbb{N}_0$):  
$$
n =  p^{k_1} + \dots + p^{k_s}, 
\qquad 
k_1 < \dots < k_s, 
\qquad 
s \le d;  
\eqno{(4)}
$$
$$
n = n_{k_1} p^{k_1} +  \dots + n_{k_s} p^{k_s},
\qquad 
k_1 < \dots < k_s, 
\qquad 
s \le d, 
\qquad 
n_{k_j} \in \{ 1, \dots, p - 1 \}. 
\eqno{(5)}
$$
Системы~(2) и (3)
являются аналогами хаосов Радемахера. 

В работе доказано, что~(2) и (3) являются системами $q$-лакунарности 
и, как следствие, системами единственности. 
Также находятся точные константы $\varepsilon$, при которых 
эти системы являются системами $\varepsilon$-единственности. 

В отличие от хаоса Радемахера, 
функции (случайные величины) из~(2) 
и~(3) являются комплекснозначными, 
а их действительные части вообще говоря не 
являются симметричными. Кроме того, 
хаос Радемахера состоит из всевозможных 
$d$-членных прозведений, составленных из 
элементов последовательности независимых функций. 
Система~(2) устроена 
подобным образом, в то время как в~(3)
элементы формируются не только из 
элементов последовательности независимых функций, 
но и из их степеней. 

В конце работы результаты об $\varepsilon$-единственности  
для систем~(2) и~(3)  
сравниваются с полученными в~[19], 
результатами подобного рода для систем из функций $\exp ( 2 \pi i n x )$, 
где $n$ берется из множеств $V^{(d)}_p$, $\widetilde{V}^{(d)}_p$ 
и $\pm \widetilde{V}^{(d)}_p$. 

\section{Основные определения, обозначения и вспомогательные факты} 

\subsection{Определения и обозначения}  

Пишем $:=$ для равенства по определению. 

Используем обозначение $a : b$ для множества $\{ a, a+1, \dots, b-1, b \}$~[20]. 

$\mu$ --- мера Лебега на множестве $[0,1)$. 

$\mathbb{N}$ означает множество натуральных, 
$\mathbb{R}$ --- действительных, 
$\mathbb{C}$ --- комплексных чисел, 
$\mathbb{N}_0 := \mathbb{N} \cup \{ 0 \}$, 
$\lfloor a \rfloor$ --- нижняя целая часть числа $a \ge 0$. 

$\mathrm{I}_E$ --- характеристическая функция (индикатор) множества $E$. 

$L_q = L_q [0,1) = 
\left\{ 
    f \colon [0,1) \to \mathbb{C} 
    \colon 
    \left( \int\limits_0^1 | f (x) |^q d \mu ) \right)^{1/q} 
    =:\| f \|_{L_q} =:\| f \|_q < \infty 
\right\}$.  

$\left\| ( c_k )_{k=1}^n \right\|_{l_2} 
:= \left( \sum\limits_{k=1}^n | c_k |^2 \right)^{1/2}$, $c_k \in \mathbb{C}$. 

Cходимость последовательностей векторов 
в пространствах $\mathbb{C}^n$ понимается 
как сходимость в евклидовой норме $\| \hspace{3pt} \|_2$ 
или в одной из эквивалентных норм.  
$\| \hspace{3pt} \|_{2,2}$ --- операторная норма матрицы, 
порожденная векторной нормой $\| \hspace{3pt} \|_2$: 
\[ 
\| A \|_{2,2} := \sup_{ \mathbf{x} \ne \mathbf{0} } 
\| A \mathbf{x} \|_2 / \| \mathbf{x} \|_2. 
\] 

До конца работы выберем и зафиксируем 
произвольное натуральное $p \ge 2$. 

$\omega := \exp ( 2\pi i / p )$.

{\it $p$-ичными интервалами} (ранга $k$) мы называем  
полуинтервалы вида 
$$
\Delta^{(k)}_m 
:= 
\big[ m p^{-k}, ( m+1 ) p^{-k} \big) \subset [0,1), 
\qquad 
k \in \mathbb{N}_0, 
\qquad 
m \in 0 : p^k - 1.  
\eqno{(6)}
$$

{\it Обобщенные функции Радемахера} на $[0,1)$ 
определяются формулой 
$$ 
R_k (x) := \omega^{ \lfloor x \cdot p^{k+1} \rfloor }, 
\qquad 
k \in  \mathbb{N}_0. 
$$
Нетрудно показать, что $R_k (x) = \omega^{ x_k }$, 
где $x_k$ --- коэффициенты из разложения 
$$  
x = \sum_{k=0}^{\infty} x_k p^{-k-1}, 
\qquad 
x_k \in 0 : p - 1,
\eqno{(7)}
$$ 
числа $x$ в бесконечную $p$-ичную дробь. 
Если таких разложений два, 
берем то, для которого $x_k=0$, 
начиная с некоторого $k$. 

Всевозможные конечные произведения функций $R_k$ 
называют {\it функциями Виленкина~---Крестенсона $VC_n$} (В--К), 
а при $p=2$ --- {\it функциями Уолша}. 
В нумерации Пэли 
$$ 
VC_n = \prod_{k=0}^{H(n)} (R_k)^{n_k}, 
$$
где $n_k$ --- коэффициенты $p$-ичного разложения числа $n \in \mathbb{N}_0$: 
$$ 
n = \sum_{k=0}^{H(n)}  n_k p^k, 
\quad 
n_k \in 0 : p - 1 
$$

При $n \le p^k - 1$ функция $VC_n$ принимает постоянное значение 
$=: VC_n ( \Delta )$ на каждом $p$-ичном интервале 
$\Delta$ ранга $k$. 
Эти значения можно собрать в 
{\it $k$-ую матрицу Виленкина~---Крестенсона}  
$\mathrm{VC}^{(k)} 
:= 
\left( VC_n ( \Delta^{(k)}_m ), \; 
n,m \in 0 : p^k - 1 \right)$  
(нумерация $\Delta^{(k)}_m$ соответствует~(6)). 
Свойства таких матриц изучались, 
напр., в~[20]. 
В частности, все $\mathrm{VC}^{(k)}$ невырожденны  
и (черта означает комплексное сопряжение)
$$
\big( 
    \mathrm{VC}^{(k)}
\big)^{-1}     
= 
\frac{1}{p^k} \, \overline{\mathrm{VC}^{(k)}}^t.  
\eqno{(8)}
$$
 

\subsection{Вспомогательные результаты}


Следующий результат доказан в~[19]. 

{\bf Лемма 1.} 
Пусть $p \ge 2$, $E_0,\dots ,E_{p-1} \subset [0,1)$, 
причем $\mu E_m \ge a$ для всех $m$. 
Тогда $\mu H \ge \frac{pa-1}{p-1}$, 
$$  
H 
:= 
\{ 
    x 
    \colon 
    \text{$x$ --- элемент минимум двух из множеств $E_0,\dots ,E_{p-1}$} 
\}.
$$  

Лемма~2 есть частный случай теоремы~6.5.1 из~[21]. 

{\bf Лемма 2.} 
Если для последовательности независимых действительнозначных 
симметрично распределенных случайных величин 
$\left \{ \xi_k \right\}_{k=1}^n$ и некоторых $q_1 \ge q_2 \ge 1$ неравенство
$$
    \left\| 
        \sum_{k=0}^n a_k \xi_k 
    \right\|_{q_1} 
    \le \kappa 
    \left\| 
        \sum_{k=0}^n a_k \xi_k 
    \right\|_{q_2}
$$
выполнено с константой $\kappa$, не зависящей от $a_k \in \mathbb{R}$, 
то с некоторой константой $C_s \ge 1$, зависящей лишь от $s$, справедливо неравенство 
$$ 
\Bigg\| 
    \sum_{0 \le k_1 < k_2 < \dots < k_s \le n} 
    a_{k_1 k_2 \dots k_s} 
    \xi_{k_1} \xi_{k_2} \dots \xi_{k_s} 
\Bigg\|_{q_1} 
\le 
C_s \kappa^s 
\Bigg\| 
\sum_{0 \le k_1 < k_2 < \dots < k_s \le n} 
    a_{k_1 k_2 \dots k_s} 
    \xi_{k_1} \xi_{k_2} \dots \xi_{k_s} 
\Bigg\|_{q_2}.
$$ 


{\bf Лемма 3.} 
Пусть функции $g_k \colon \{0, \dots, p-1\} \rightarrow \mathbb{R}$,  
$k \in \mathbb{N}_0$, заданы формулой $g_k (m) = \alpha_k^m$.  
Тогда функции $f_k(x) := g_k(x_k)$ 
образуют последовательность независимых функций. Здесь 
$x_k$ --- коэффициенты из разложения~(7). 

{\bf Доказательство.} 
Достаточно доказать для всех 
$n \in \mathbb{N}$ и $\varepsilon_k \in \mathbb{R}$ 
равенство 
$$ 
\mu 
\{ 
    x \in [0,1) \colon f_0 (x) = \varepsilon_0, \, \dots, \, f_n (x) = \varepsilon_n 
\}  
= 
\prod\limits_{k=0}^n 
\mu 
\{ 
    x \in [0,1) \colon f_k (x) = \varepsilon_k 
\}. 
\eqno{(9)}
$$ 
Пусть $\varepsilon_k$ встречается $M(k)$ раз среди $\alpha_k^m$. 
Тогда множество $\{ x \in [0,1) \colon f_k (x) = \varepsilon_k \}$ 
является дизъюнктным объединением $M(k) p^k$ 
$p$-ичных интервалов ранга $k+1$ и его мера есть $M(k) p^k / p^{k+1} = 
M(k) / p$.  
Тогда вся правая часть~(9) равна $p^{-(n+1)} \prod\limits_{k=0}^n M(k)$. 

Вычислим левую часть~(9). 
Если $b_k \in 0 : p-1$, то множество 
тех $x$, коэффициенты $x_k$ из разложения~(7) 
которых совпадают с $b_k$ для всех $k=0,\dots ,n$ 
образуют некоторый $p$-ичный интервал ранга $p^{-(n+1)}$. 
Тогда множество из левой части~(9) 
является дизъюнктным объединением $\prod\limits_{k=0}^n M(k)$ 
таких $p$-ичных интервалов и его мера  
$p^{-(n+1)} \prod\limits_{k=0}^n M(k)$ совпадает 
с правой частью~(9). 
Лемма доказана.

\section{Основные результаты} 

Сначала покажем, что системы~(2) 
и~(3) являются $q$-лакунарными и, 
как следствие, системами единственности. 
    
{\bf Теорема 1.} 
Для каждого $q > 2$ система функций~(3) является 
системой $q$-лакунарности, т.е. для всех натуральных $N$  
$$
\Bigg\| 
    \sum_{ n \in \widetilde{V}^{(d)}_p \cap [ 1, N ] } 
    c_n VC_n 
\Bigg\|_q 
\le 
C  
\big\| 
    ( c_n )_{ n \in \widetilde{V}^{(d)}_p \cap [ 1, N ] }
\big\|_{l_2}, 
\eqno{(10)}
$$ 
$C = C ( p, d, q )$ не зависит от $c_n \in \mathbb{C}$. 
Как следствие, (2) также 
есть система $q$-лакунарности. 

{\bf Доказательство.}
Положим $\widetilde{W}^{(s)}_p := \widetilde{V}^{(s)}_p \setminus 
\widetilde{V}^{(s-1)}_p$, т.е. $\widetilde{W}^{(s)}_p$ 
состоит из всех $n \in \mathbb{N}$, 
у которых ровно $s$ ненулевых $p$-ичных коэффициентов. 
Очевидно, 
$$ 
\widetilde{V}^{(d)}_p 
= 
\bigsqcup\limits_{s=1}^d 
\widetilde{W}^{(s)}_p.
$$ 

Возьмем и зафиксируем любое натуральное $L$, 
а вслед за ним произвольный набор чисел $j_0,\dots ,j_L$, взятых из множества 
$\{ 1, \dots, p-1 \}$. 
Рассмотрим множество 
$$
A^{(s)}_{ j_0 \dots j_L } 
:= 
\big\{ 
    n \in \widetilde{W}^{(s)}_p 
    \colon  
    n = j_{k_1} p^{k_1} + \dots + j_{k_s} p^{k_s}
\big\}, 
$$ 
состоящее из тех натуральных $n$, у которых 
ровно $s$ ненулевых $p$-ичных коэффициентов, 
причем каждый ненулевой $n_k$ совпадает с $j_k$. 
Очевидно, $n < p^{L+1}$ для таких $n$.  
 
Случайные величины $\mathrm{Re} \, R_k^{j_k}$ 
и $\mathrm{Im} \, R_k^{j_k}$ действительнозначны. 
Они симметричны в случае четного основания $p$ и  
как минимум вторая --- в случае нечетного. 
В любом случае $\mathrm{Re} \, R_k^{j_k}$ 
можно представить, как сумму $p-1$ симметрично распределенной случайной величины: 
$$
\mathrm{Re} \, R_k^{j_k} (x)
= 
\sum_{m=0}^{p-1} 
\cos \left( \frac{ 2 \pi m j_k }{p} \right) 
\mathrm{I}_{ \left[ m/p, (m+1)/p \right) } (p^k x) 
$$ 
$$ 
= 
\sum_{m=1}^{p-1} 
\left( 
    \cos \left( \frac{ 2 \pi m j_k }{p} \right) 
    \mathrm{I}_{ \left[ m/p, (m+1)/p \right) } (p^k x)  
    - 
    \cos \left( \frac{ 2 \pi m j_k }{p} \right) 
    \mathrm{I}_{ \left[ 0, 1/p \right] } (p^k x) 
\right),  
\eqno{(11)}
$$  
нужное разложение содержится в нижней строке~(11). Здесь 
индикаторы $\mathrm{I}_{ \left[ m/p, (m+1)/p \right) }$ $1$-периодично 
продолжены на правую полуось. 

Для каждого $k=0,\dots ,L$ обозначим $f_{k, \alpha }(x)$ 
функции $\mathrm{Re} \, R_k^{j_k} (x)$ и $\mathrm{Im} \, R_k^{j_k} (x)$, 
если $p$ четно (тогда $\alpha$ пробегает два значения), 
и функцию $\mathrm{Im} \, R_k^{j_k} (x)$, а также 
функции в скобках в~(11), 
если $p$ нечетно (тогда $\alpha$ пробегает $p$ значений). 
Кроме того, для комплексного $c_n$ 
обозначим $a_{n, \alpha}$ его действительную либо мнимую часть. 

Неравенства Минковского и Коши~---Буняковского дают 
$$
\Bigg\| 
    \sum_{ n \in A^{(s)}_{j_0 \dots j_L } } 
    c_n VC_{n} 
\Bigg\|_q^2 
\le 
\Bigg( 
    \sum_{ \alpha_0, \alpha_1, \dots, \alpha_s } 
    \Bigg\| 
        \sum_{ 0 \le k_1 < \dots < k_s \le L } 
        a_{k_1, \dots, k_s, \alpha_0} 
        f_{k_1, \alpha_1} \cdot \dots \cdot f_{k_s, \alpha_s} 
    \Bigg\|_q 
\Bigg)^2
$$ 
$$ 
\le 
p^s 
\sum_{ \alpha_0, \alpha_1, \dots, \alpha_s } 
\Bigg\| 
    \sum_{ 0 \le k_1 < \dots < k_s \le L } 
    a_{k_1, \dots, k_s, \alpha_0} 
    f_{k_1, \alpha_1} \cdot \dots \cdot f_{k_s, \alpha_s} 
\Bigg\|_q^2,   
\eqno{(12)}
$$ 
где $0 \le k_1 < \dots < k_s \le L$, $n= j_{k_1} p^{k_1} + \dots + j_{k_s}p^{k_s}$, 
$a_{k_1, \dots, k_s, \alpha_0} = a_{n, \alpha_0}$. 
По лемме~3 случайные величины $f_{ k, \alpha }$   
с различными $k$ действительнозначны, ограничены, независимы 
и имеют нулевое математическое ожидание.  
Поэтому~[22,\,23]
для любого набора $\beta_0,\dots ,\beta_L$ 
справедливо следующее: 
$$ 
\Bigg\| 
    \sum_{ k = 0 }^n  
    a_k f_{ k, \beta_k } 
\Bigg\|_q 
\le 
\kappa
\Bigg\| 
    \sum_{ k = 0 }^n  
    a_k f_{ k, \beta_k } 
\Bigg\|_2, 
\quad 
\text{$\kappa = \kappa (q)$ не зависит от $a_k \in \mathbb{R}$}.
$$
Т.к. случайные величины $f_{ k, \alpha }$ еще и симметричны, 
из неравенства выше и леммы~2 получаем, что 
правая часть~(12) не превосходит 
$$
C_1 
\sum_{ \alpha_0, \alpha_1, \dots, \alpha_s } 
\Bigg\| 
    \sum_{ 0 \le k_1 < \dots < k_s \le L } 
    a_{k_1, \dots, k_s, \alpha_0} 
    f_{k_1, \alpha_1} \cdot \dots \cdot f_{k_s, \alpha_s} 
\Bigg\|_2^2, 
\qquad 
C_1 = C_1 ( p, d, q ). 
\eqno{(13)}
$$
Т.к. случайные величины $f_{ k, \alpha }$ симметричны 
и независимы, они некоррелированы. 
Их $L_2$-нормы не превосходят $1$. 
Следовательно, выражение в~(13) 
не превосходит 
$$ 
C_1  
\sum_{ \alpha_0, \alpha_1, \dots, \alpha_s } 
\sum_{ 0 \le k_1 < \dots < k_s \le L } 
| a_{k_1, \dots, k_s, \alpha_0} |^2 
\le 
C_1 p^{s+1} 
\sum_{ 0 \le k_1 < \dots < k_s \le L } 
| a_{k_1, \dots, k_s, \alpha_0} |^2 
$$
$$ 
\le 
C_2 
\sum_{ n \in A^{(s)}_{j_0 \dots j_L } } 
| c_n |^2, 
\qquad 
C_2 = C_2 ( p, d, q ).  
\eqno{(14)}
$$
Из~(12)--(14) получаем 
$L_2$-$L_q$-неравенство Хинчина для систем $A^{(s)}_{ j_0 \dots j_L } \subset \widetilde{W}^{(s)}_p$: 
$$
\Bigg\| 
    \sum_{ n \in A^{(s)}_{j_0 \dots j_L} \cap [1, N] } 
    c_n VC_n 
\Bigg\|_q 
\le 
C_2 
\sum_{ n \in A^{(s)}_{j_0 \dots j_L } } 
| c_n |^2. 
\eqno{(15)}
$$
Распространим~(15) на систему $\widetilde{W}^{(s)}_p$. 
Запишем равенство 
$$ 
\sum_{ n \in \widetilde{W}^{(s)}_p \cap [1, N] } 
c_n VC_n 
= 
\frac{1}{ (p-1)^{L+1-s} } 
\sum_{j_0, \dots, j_L } 
\sum_{ n \in A^{(s) }_{ j_0 \dots j_L } \cap [1, N] } 
c_n VC_n, 
\qquad 
p^L \le N < p^{L+1}.  
\eqno{(16)}
$$
В~(16) сумма слева разбита  
на $(p-1)^{L+1}$ сумм, 
в каждую из которых входят функции В--К $VC_n$ 
вида $(R_{k_1})^{j_{k_1}} \dots (R_{k_s})^{j_{k_s}}$ 
с фиксированным $s$, а степени зависят только от номера $k$ функции $R_k$. 
Множитель $(p-1)^{s-L-1}$ возникает из-за того, что каждая такая функция входит 
в $(p-1)^{L+1-s}$ разных сумм 
(другими словами, каждое $n \in \widetilde{W}^{(s)}_p \cap [1, N]$ 
попадает в в $(p-1)^{L+1-s}$ множеств вида $A^{(s) }_{ j_0 \dots j_L } \cap [1, N]$). 
Применяя~(15)--(16) и неравенство Минковского, 
а также используя ортогональность системы В--К, 
получаем:  
$$ 
\Bigg\| 
    \sum_{ n \in \widetilde{W}^{(s)}_p \cap [1, N] } 
    c_n VC_n 
\Bigg\|_q^2 
\le 
\frac{1}{ (p-1)^{ 2 ( L+1-s ) } } 
\Bigg( 
    \sum_{ j_0, \dots, j_L } 
    \bigg\| 
        \sum_{ n \in A^{(s)}_{ j_0 \dots j_L }  \cap [1, N] } 
        c_n VC_n 
    \bigg\|_q 
\Bigg)^2 
$$
$$ 
\le 
\frac{ (p-1)^{ L+1 } }{ (p-1)^{ 2 ( L+1-s ) } } 
\sum_{ j_0, \dots, j_L } 
\bigg\| 
    \sum_{ n \in A^{(s)}_{ j_0 \dots j_L }  \cap [1, N] } 
    c_n VC_n 
\bigg\|_q^2 
$$
$$ 
\le 
\frac{ (C_2)^2 }{ (p-1)^{L+1-2s} } 
\sum_{ j_0, \dots, j_L } 
\sum_{ n \in A^{(s)}_{j_0 \dots j_L } } 
| c_n |^2
=  
\frac{ (C_2)^2 (p-1)^{L+1-s} }{ (p-1)^{L+1-2s} } 
\sum_{ n \in \widetilde{W}^{(s)}_p  \cap [1, N] } 
| c_n |^2 
$$
$$ 
= 
(C_2)^2 (p-1)^{s} 
\sum_{ n \in \widetilde{W}^{(s)}_p  \cap [1, N] } 
| c_n |^2.  
$$
Наконец, с учетом цепочки выше, 
для системы $\widetilde{V}^{(d)}_p = \bigsqcup\limits_{s=1}^{d} 
\widetilde{W}^{(s)}_p$ верно следующее:  
$$ 
\Bigg\| 
    \sum_{ n \in \widetilde{V}^{(d)}_p \cap [1, N] }  
    c_n VC_n 
\Bigg\|_q^2 
\le 
\Bigg( 
    \sum_{s=1}^{d} 
    \Bigg\| 
        \sum_{ n \in \widetilde{W}^{(s)}_p \cap [1, N] } 
         c_n VC_n  
    \Bigg\|_q  
\Bigg)^2 
\le 
d  
\sum_{s=1}^{d} 
\Bigg\| 
    \sum_{ n \in \widetilde{W}^{(s)}_p \cap [1, N] } 
    c_n VC_n  
\Bigg\|_q^2  
$$
$$ 
\le 
d  
(C_2)^2 (p-1)^d  
\sum_{s=1}^{d} 
\sum_{ n \in \widetilde{W}^{(s)}_p \cap [1, N] } 
| c_n |^2 
= 
C^2   
\sum_{ n \in \widetilde{V}^{(d)}_p \cap [1, N] }  
| c_n |^2, 
\qquad 
C = C ( p, d, q ), 
$$
а это равносильно~(10). 

Далее мы находим точные константы для $\varepsilon$-единственности 
систем~(2) и~(3). 
Следующие леммы~4 и 5 обобщают 
в разных направлениях обобщают теорему~3 из~[12]. 
Содержащиеся в них утверждения служат 
базой индукции для теорем~2 и 3. 
Техника доказательства лемм~4--5 и теорем~2--5 
базируется на идеях из~[1,\,12,\,19].


{\bf Лемма 4.} 
Если ряд по системе обобщенных функций Радемахера 
$$ 
\sum\limits_{k=0}^{\infty} c_k R_k(x)
\eqno{(17)}
$$
сходится к постоянной $C$ на множестве $E \subset [0,1)$ меры 
$\mu(E) > \frac{1}{p}$, то $C=0$ и все $c_{k}=0.$

{\bf Доказательство.}
Во-первых, если ряд~(17) 
сходится на множестве $E$ положительной меры $E$, то он сходится п.в. на $[0,1)$. 
Доказательство этого факта практически дословно 
повторяет доказательство теоремы 1 в~[12], 
нужно лишь степени двойки поменять на степени $p$. 

Дальнейшее доказательство проведем от противного. 
Предположим, что существует не тождественно нулевой ряд~(17), 
сходящийся к нулю на множестве $E$ меры $\mu(E) > p^{-1}$. 
Возьмем первый отличный от нуля коэффициент $c_{\widetilde{k}}$. 
Для всех 
$m_1$, $m_2 \in 0 : p-1$, $m_1 \ne m_2$, 
рассмотрим выражение 
$$
d_{m_1, m_2} (x) 
:= 
\sum_{k=0}^\infty 
c_k R_k 
\big( 
    x + {m_1}{ p^{ - \widetilde{k} - 1 } } 
\big) 
- 
\sum_{k=0}^\infty 
c_k R_k 
\big( 
    x + {m_2}{ p^{ - \widetilde{k} - 1 } } 
\big). 
\eqno{(18)}
$$ 
При $k > \widetilde{k}$ члены обоих рядов из~(18)  
совпадают между собой в силу $p^{-k}$-периодичности функции 
$R_k$. Отсюда и из выбора 
$c_{\widetilde{k}}$ 
вытекает равенство 
$$
d_{m_1, m_2}(x) 
= 
c_{\widetilde{k}} 
\big[ 
    R_{\widetilde{k}} 
    \big( 
        x + {m_1}{ p^{ -\widetilde{k} - 1} } 
    \big)
    - 
    R_{\widetilde{k}} 
    \big( 
        x + {m_2}{ p^{ -\widetilde{k} - 1} } 
    \big) 
\big].
\eqno{(19)}
$$
Из определения обобщенных функций Радемахера вытекает, 
что выражение в квадратных скобках в~(19)  
не равно нулю. 
Поэтому $d_{m_1, m_2}(x) \ne 0$ для п.в. $x \in [0, 1)$.

С другой стороны, для всех $m \in 0 : p - 1$ 
сумма ряда 
$\sum_{k=0}^\infty c_k R_k 
\big( x + {m}{ p^{ -\widetilde{k} - 1 } } \big)$ 
равна $C$ на множестве 
$$
E_m := E - {m}{ p^{ - \widetilde{k} - 1 } } 
\, ( \mathrm{mod} \, 1 ) 
\eqno{(20)}
$$ 
меры большей $1 / p$. 
Если предположить, что $\mu( E_{m_1} \cap E_{m_2} ) = 0$ 
при всех $m_1 \ne m_2$, то 
$$
\mu \left( \bigcup\limits_{m=0}^{p-1} E_m \right) 
= 
\sum\limits_{m=0}^{p-1} \mu( E_m ). 
\eqno{(21)}
$$ 
Но левая часть~(21) не больше $1$, 
а правая больше, т.к. каждое из $p$ 
слагаемых больше $p^{-1}$. 
Значит, 
$\mu ( E_{m_1} \cap E_{m_2} ) \ne 0$ 
для некоторых $m_1 \ne m_2$.  
Следовательно,  
$d_{m_1, m_2}(x) = 0$ 
на множестве 
$E_{m_1} \cap E_{m_2}$ положительной меры. 
Это противоречит тому, 
что $d_{m_1, m_2}(x) \ne 0$ п.в. на $[0, 1)$. 
Полученное противоречие доказывает теорему.

{\bf Лемма 5.} 
Если $\mu(E) > 1- \frac{1}{p}$, $E \subset [0,1)$, и ряд 
$$
\sum\limits_{n \in \widetilde{V}^{(1)}_p } c_{n} VC_{n}(x)
\eqno{(22)}
$$
сходится к некоторой постоянной $C$ на множестве $E$, 
то $C=0$ и все $c_n=0.$

{\bf Доказательство.}
Предположим, что существует ненулевой ряд~(22), 
сходящийся к нулю на множестве $E$ меры $\mu(E) > \frac{p-1}{p}$. 
Из всех $n \ne 0$ с $c_n \neq 0$ выберем $n$ с минимальным 
$k_1$ из формулы~(5). 
Обозначим соответственно $\widetilde{n}$ и $\widetilde{k}$ 
такие $n$ и $k_1$. 

Mера множества   
$$
X 
= 
\bigcap\limits_{m = 0}^{p-1}
E_m, 
\eqno{(23)}
$$
где $E_m$ введено в~(20), положительна. 
Фиксируем произвольное $x \in X$. 
Т.к. всякое $n \in \widetilde{V}_p^{(1)}$ имеет вид ${j p^k}$, 
где $j \in 1 : p - 1$, 
имеем $\widetilde{n} = j p^{\widetilde{k}}$. 
Учитывая этот факт,  
из сходимости ряда~(22) к $C$ на множестве $E$ 
получаем:   
$$ 
\sum_{j=1}^{p-1} 
c_{ j p^{\widetilde{k}} } R_{ \widetilde {k} }^j 
( x + m p^{ - \widetilde{k} - 1 } ) 
+ 
\sum_{ n \in A_0 } 
c_n VC_n ( x + m p^{ - \widetilde{k} - 1 } ) = C, 
\qquad 
m \in 0 : p - 1, 
\eqno{(24)}
$$
$$
A_0 := 
\big\{ 
    n \in \widetilde{V}^{(1)}_p 
    \colon n = 0 \ 
    (\bmod \ p^{\widetilde{k}+1}) 
\big\}.
$$

Для всех $m$ и $j \in 0 : p-1$ 
из определения обобщенных функций Радемахера 
вытекает следующее: 
$$
R_{ \widetilde{k}}^j ( x + m p^{ - \widetilde{k} - 1 } ) 
= 
R_{ \widetilde{k} }^j (x) \, \omega^{jm} 
\quad \text{и} \quad 
VC_n(x) = VC_n ( x + mp^{-\widetilde{k}-1} ), \; n \in A_{0}.
$$
Тогда для всех $x \in E$ равенство из~(24) можно записать 
так: 
$$
\lim\limits_{ N \rightarrow \infty } 
\mathrm{VC}^{(1)} \textbf{S}_N (x) 
= 
\textbf{C}. 
\eqno{(25)}
$$ 
Здесь $\mathrm{VC}^{(1)}$ --- матрица Виленкина--Крестенсона 
с номером $1$, а $S_N$ и $C$ --- $(p\times 1)$-матрицы 
$$ 
\textbf{S}_N (x)
:= 
\bigg( 
    \sum_{ n \in A_{0} \colon n \le N } 
    c_n VC_n(x)  
    \;\;\;
    c_{ p^{\widetilde{k}} } R_{ \widetilde{k} } (x)  
    \;\;
    \dots   
    \;\;  
    c_{ (p-1) p^{\widetilde{k}} } 
    R^{p-1}_{ \widetilde{k} } (x) 
\bigg)^t, 
$$
$\textbf{C} := (C \; C \; \dots \; C)^t$. 
Обозначим $\textbf{C}_N (x) := \mathrm{VC}^{(1)} \textbf{S}_N (x)$.  
Тогда $\textbf{S}_N (x) = \big( \mathrm{VC}^{(1)} \big)^{-1} 
\textbf{C}_N (x)$. Т.к., согласно~(25),  
$\lim\limits_{ N \rightarrow \infty } 
\big\| \textbf{C}_N (x) - \textbf{C} \big\|_2 = 0$ при $x \in E$, 
для таких $x$ 
$$
\big\| 
    \textbf{S}_N (x) 
    - 
    \big( \mathrm{VC}^{(1)} \big)^{-1} \textbf{C} 
\big\|_2 
\le 
\big\| 
    \big( \mathrm{VC}^{(1)} \big)^{-1} 
\big\|_{2,2} 
\big\|
    \textbf{C}_N (x) - \textbf{C} 
\big\|_2 \rightarrow 0,   
\;\; \rightarrow \infty. 
$$  
С учетом~(8) получаем: 
$$ 
\lim\limits_{ N \to \infty} 
\textbf{S}_N (x) 
= 
\big( \mathrm{VC}^{(1)} \big)^{-1} \textbf{C} 
= 
\frac{1}{p} \, \overline {\mathrm{VC}^{(1)} }^t \textbf{C} 
= 
(C \; 0 \; \dots \; 0)^t. 
$$
Тем самым для каждого $j \in 1 : p-1$ и $x \in E$ 
верно равенство $c_{ j p^{\widetilde{k}} } R_{ \widetilde{k} }^j (x) = 0$, 
откуда $c_{ j p^{\widetilde{k}} } = 0$. 
Это противоречит выбору $\widetilde{k}$. 
Полученное противоречие доказывает теорему.


{\bf Теорема 2.} 
Если на множестве $E \subset [0,1)$ c $\mu E > 1 - \left(\frac{p-1}{p} \right)^{d}$ 
ряд 
$$
\sum_{n \in V^{(d)}_p } c_{n} VC_{n}(x)
\eqno{(26)}
$$
сходится к некоторой постоянной $C$, то $C=0$ и все $c_{n}=0.$

{\bf Доказательство.}
Будем доказывать индукцией по $d$. При $d = 1$ утверждение верно 
согласно лемме~4. 
Допустим, оно верно для $d-1$ и не верно для $d$ 
и существует не тождественно нулевой ряд~(26), 
который сходится к $C$ на множестве $E$ меры 
$\mu E > 1 - \left(\frac{p-1}{p}\right) ^{d}$. 
Выберем $\widetilde{n}$ и $\widetilde{k}$ так, 
как в лемме~5, 
но используя разложение~(4) вместо~(5). 
Пусть $F_{mj} = E_m \cap E_j$, где 
множества справа определены в~(20). 
Если $X$ --- объединение всевозможных $F_{ij}$, 
то (лемма~1)  
$$ 
\mu ( X ) 
> 
\frac{ p \mu(E) - 1 }{p-1} 
> 
\frac{ 1 }{ p - 1 } 
\left(p \left( 1 - \left(\frac{p-1}{p} \right)^{d} \right) -  1 \right) 
= 1 - \left(\frac{p-1}{p}\right)^{d-1}.
$$

Для произвольных $F_{mj}$ и $x \in F_{mj}$ рассмотрим выражение 
$$
\sum_{ n \in V^{(d)}_p } 
c_n VC_n 
\big( x + m{p^{-\widetilde{k}-1}} \big) 
- 
\sum_{ n \in V^{(d)}_p } 
c_n VC_n
\big( x + j{p^{-\widetilde{k}-1}} \big)
\eqno{(27)}
$$ 
Оно равно нулю, так как обе суммы равны $C$. 
С другой стороны, оно не изменится, 
если из обеих сумм~(27) одновременно удалить члены 
с одинаковыми номерами $n$ такими, 
что в $p$-ичном разложении~(4) числа 
$n$ либо $k_1 > \widetilde{k}$, либо $k_1 < \widetilde{k}$   
(в первом случае мы удаляем разность одинаковых, 
во втором --- нулевых чисел). Тогда 
$$
\sum\nolimits^{\prime}  
c_n 
\big[ 
    VC_n 
    \big( x + m { p^{-\widetilde{k}-1} } \big) 
    -
    VC_n 
    \big( x + j { p^{-\widetilde{k}-1} } \big)
\big] 
= 0,  
\eqno{(28)}
$$   
где $\sum^{\prime}$ распространяется на 
$$ 
n \in V^{(d)}_p 
\colon \,  
n - p^{\widetilde{k}} \in V^{(d-1)}_p 
\;\; \wedge \;\;  
n = p^{\widetilde{k}} + a(n) p^{ \widetilde{k} + 1 }, 
\;\; a(n) \in \mathbb{N}_0. 
$$ 

Преобразуем выражение в квадратных скобках в~(28): 
$$ 
VC_n 
\big( x + m { p^{-\widetilde{k}-1} } \big) 
-
VC_n 
\big( x + j { p^{-\widetilde{k}-1} } \big)
$$ 
$$ 
= 
VC_{ p^{\widetilde{k}} } 
\big( x + m { p^{-\widetilde{k}-1} } \big) 
VC_{ a(n) p^{ \widetilde{k} + 1 } }
\big( x + m { p^{-\widetilde{k}-1} } \big) 
- 
VC_{ p^{\widetilde{k}} } 
\big( x + j { p^{-\widetilde{k}-1} } \big) 
VC_{ a(n) p^{ \widetilde{k} + 1 } }
\big( x + j { p^{-\widetilde{k}-1} } \big) 
$$ 
$$ 
= 
R_{ \widetilde{k} } 
\big( x + m { p^{-\widetilde{k}-1} } \big) 
VC_{ a(n) p^{ \widetilde{k} + 1 } }
( x ) 
- 
R_{ \widetilde{k} } 
\big( x + j { p^{-\widetilde{k}-1} } \big) 
VC_{ a(n) p^{ \widetilde{k} + 1 } }
( x )  
$$
$$ 
= 
R_{ \widetilde{k} } 
( x ) 
VC_{ a(n) p^{ \widetilde{k} + 1 } }
( x ) 
( \omega^m - \omega^j ).  
$$
Т.к. $\omega^m - \omega^j \ne 0$, цепочки выше и~(28) дают 
$$
\sum\nolimits^{\prime}  
c_n 
VC_{ a(n) p^{ \widetilde{k} + 1 } } (x) = 0.
\eqno{(29)}
$$
Ряд из~(29) есть  
не тождественно нулевой ряд вида 
$\sum\nolimits_{n \in V^{(d-1)}_p } d_{n} VC_{n}(x)$, 
один и тот же для всех $F_{ij}$, 
который сходится к нулю на множестве $X$ меры, 
большей $1 - \left(\frac{p-1}{p}\right) ^{d-1}$. 
Получаем противоречие с предположением индукции. 
Теорема доказана.


{\bf Теорема 3.} 
Если на множестве $E \subset [0,1)$ 
c $\mu E > 1 - p^{-d}$ ряд
$$
\sum_{n \in \widetilde{V}_p^{(d)}} c_n VC_n(x)
\eqno{(30)}
$$
сходится к некоторой постоянной $C$, то $C=0$ и все $c_{n}=0.$

{\bf Доказательство.}
Доказательство проведем индукцией по $d$.
При $d=1$ утверждение следует из леммы~5. 
Допустим, утверждение теоремы верно для $d-1$; 
докажем, что оно верно и для $d$.  
Предположим, что это не так 
и существует ряд~(30), 
не все коэффициенты которого нулевые, 
сходящийся к $C$ на множестве $E$ меры $\mu E > 1 - p^{-d}$. 
Выберем $\widetilde{n}$ и $\widetilde{k}$ так, 
как в теореме~3.  
Построим множество $X$ с помощью формулы~(23); 
его мера больше, чем $1 - p \cdot p^{-d} = 1 - p^{-(d-1)}$. 
Фиксируем произвольное $x \in X$. 
Учитывая выбор $\widetilde{n}$ и $\widetilde{k}$,  
из сходимости ряда~(30) к $C$ на множестве $E$ 
получаем:   
$$  
\sum_{ n \in A_0 \sqcup \dots \sqcup A_{p-1} } 
c_n VC_n 
\big( 
    x + m p^{ - \widetilde{k} - 1 } 
\big) = C, 
\qquad 
m = 0 : p-1,
\eqno{(31)}
$$ 
$$ 
A_j  
:= 
\big\{ 
    n \in \widetilde{V}^{(d)}_p 
    \colon 
    n - j p^{\widetilde{k}} = 0 \; ( \mathrm{mod} \, p^{\widetilde{k}+1} )  
\big\}.  
$$

Из определения функций В--К вытекает следующее: 
$$ 
VC_n 
\big( 
    x + m p^{-\widetilde{k}-1} 
\big) 
= 
\omega^{m j } VC_{ n - j p^{\widetilde{k}} } (x) VC_{ j p^{\widetilde{k}} }(x), 
\quad 
m, \, j \in 0 : p-1, 
\quad 
n \in A_j. 
$$ 
Тогда~(31) записывается в виде~(25), 
где $\textbf{S}_N$ и $\textbf{C}$ --- $(p\times 1)$-матрицы,  
$$
\textbf{S}_N  
= 
\bigg( 
    \sum_{n \in A_j \colon n \le N } c_n VC_{ j p^{\widetilde{k}} } (x) VC_{ n -j p^{\widetilde{k}} } (x) 
\bigg)^t_{ j \in 0 : p-1}, 
\qquad 
\textbf{C} = (C \; C \; \dots \; C)^t.
$$
Рассуждая, как при доказательстве леммы~5, 
получаем, что для каждого $j=0 : p-1$ ряд
$$
\sum_{n \in A_j \colon n \leqslant N} 
c_n VC_{ j p^{\widetilde{k}} } (x) 
VC_{ n -j p^{\widetilde{k}} } (x)
$$
сходится к $j$-ой компоненте вектора $W^{-1} \textbf{C}$. 
Так как 
$$ 
( C_0 \; C_1 \; \dots \; C_{p-1} )^t 
= 
\frac{1}{p} \overline{W}^t ( C \; C \; \dots \; C)^t 
= 
(C \; 0 \; \dots \; 0),
$$ 
для всех $x \in X$ ряд 
$$
\sum_{ n \in A_j } 
c_n VC_{ n - j p^{\widetilde{k}} } (x)
\eqno{(32)}
$$
сходится к нулю для всех $j=1 : p-1$. 
Заметим, что $n - j p^{\widetilde{k}} \in \widetilde{V}_p^{(d-1)}$ и 
$\mu(X) > 1 - p^{-d+1}$. По предположению индукции все ряды~(32)тождественно равны 0. Это противоречит тому, что $c_n \ne 0$ для некоторого 
$j \in 1 : p-1$ и $n \in A_j$, 
согласно построению множеств $A_j$. 
Противоречие завершает индукционный 
переход и доказывает теорему.  

Из теорем~2 и 3 вытекает 

{\bf Следствие 1.} 
При $\mu(E) < \left( \frac{p-1}{p} \right)^{d}$ 
множество $E\subset [0,1)$ является {\it множеством единственности} 
для системы~(2), 
а при $\mu(E) < p^{-d}$ --- для системы~(3). 
Это означает, что сходящийся вне $E$ к нулю ряд 
по соответствующей системе содержит лишь нулевые коэффициенты. 

{\bf Следствие 2.} 
(2) --- система $\varepsilon$-единственности при 
$\varepsilon \le \left( \frac{p-1}{p} \right)^{d}$, 
а~(3) --- при 
$\varepsilon \le  p^{-d}$.

Неравенства в условиях на меру в 
теоремах~2 и 3 
и в следствии~1, 
а также на $\varepsilon$ в следствии~2 являются точными, 
т.к. многочлены  
$$ 
\prod_{k=0}^{d-1} ( 1 - R_k(x) ) - 1 
\;\; \text{и} \;\; 
\prod_{k=0}^{d-1} ( 1 + R_k(x) + \dots + R_k^{p-1}(x) ) - 1
$$ 
обращаются в ноль на множествах меры $1 - \left( \frac{p-1}{p} \right)^{d}$  
и $1 - p^{-d}$, 
соответственно. 

Сравним полученные результаты с аналогичными, 
где вместо функций Виленкина--Крестенсона 
берутся тригонометрические функции. 
Так, в~[19]
показано, что 
{\it 
система функций $\{ \exp ( 2 \pi i n x ), \, n \in V_p^{(d)} \}$
является системой $\varepsilon$-единственности 
при $\varepsilon \le \left( \frac{p-1}{p} \right)^{d-1}$} 
(о точности значения $\varepsilon$ мы не знаем). 
Этот результат дает более широкий класс 
систем $\varepsilon$-единственности, 
чем для системы~(2) 
(следствие~2). 
В тоже время системы
$\{ \exp ( 2 \pi i n x ), \, n \in \widetilde V_p^{(d)} \}$ 
и 
$\{ \exp ( 2 \pi i n x ), \, n \in \pm \widetilde V_p^{(d)} \}$
--- системы единственности 
при 
$\varepsilon \le  p^{-(d-1)}$ 
и 
$\varepsilon \le  p^{-2(d-1)}$, 
соответственно 
(и снова мы не знаем о точности значения $\varepsilon$), 
причем для первой системы результат дает более широкий, 
для второй --- более узкий класс 
систем единственности, 
чем для подсистемы~(3) 
системы В--К (следствие~2).

\end{document}